\documentclass[11pt,twoside]{article}

\usepackage{amsmath,amsfonts,amsthm,graphicx,amscd,amssymb,latexsym,mathrsfs,graphics,color,longtable,multicol,verbatim,lscape}
\usepackage[mathscr]{euscript}
\usepackage{enumerate}

\usepackage[all]{xy}

\usepackage[bookmarks,colorlinks,linkcolor=black,citecolor=black]{hyperref}

\numberwithin{equation}{section}

\font\ita=cmti9

\font\ita=cmti9   

\newtheorem{thm}{Theorem}[subsection]
\newtheorem{prop}[thm]{Proposition}
\newtheorem{cor}[thm]{Corollary}
\theoremstyle{definition}
\newtheorem{defn}[thm]{Definition}
\newtheorem{rem}[thm]{Remark}

\title{{\bf Cohomology of some families of Lie algebras and quadratic Lie algebras}
    \footnotetext{
	 {\bf MSC (2010):} \footnotesize{Primary 17B56, 17B60, 16E40; Secondary 17B30, 16W25.\\
 	 {\bf Financial Support:} The paper was supported by the University of Economics and Law, VNU-HCMC.
	}}} 
\author{{\bf Cao Tran Tu Hai$^*$, Duong Minh Thanh$^\dagger$ , 
Le Anh Vu$^{**}$}\\
	\vspace{-0.2cm}{\ita $^*$  Le Quy Don High School for the Gifted, Ninh Thuan, Vietnam}\\
	\vspace{-0.2cm}{\ita Email: tuhai.thptlequydon@ninhthuan.edu.vn}\\
	\vspace{-0.2cm}{\ita $^\dagger$ Ho Chi Minh city University of Education, Vietnam}\\
	\vspace{-0.2cm}{\ita Email: thanhdmi@hcmue.edu.vn}\\
	\vspace{-0.2cm}{\ita ${}^{**}$University of Economics and Law, VNU-HCMC, Vietnam}\\
	\vspace{-0.2cm}{\ita Email: vula@uel.edu.vn}}

\begin{document}
\date{}
\maketitle
\label{firstpage}

\begin{abstract}
The paper studies the cohomology of Lie algebras and quadratic Lie algebras. Firstly, we propose to describe the cohomology of  $MD(n,1)$-class which was introduced in \cite{LHNCN16}. This class contains Heisenberg Lie algebras. In 1983, L. J. Santharoubane \cite{San83} computed the cohomology of Heisenberg Lie algebras. In this paper, we will completely describe the cohomology of the other ones of $MD(n, 1)$-class. Finally, we will be concerned about the cohomology of quadratic Lie algebras. In 1985, A. Medina and P. Revoy \cite{MR85} computed the second Betti number of the generalized real diamond Lie algebras. We will compute in this paper the second Betti number of the generalized complex diamond Lie algebras by using the super-Poisson bracket.
\end{abstract}

\markboth{Cohomology of Lie algebras}{Cohomology of quadratic Lie algebras}
\section*{Introduction}
	Cohomology of Lie algebra is a special case of cohomology of algebras which is inspired by de Rham cohomology on the corresponding Lie group and has both algebraic and geometric flavor. We restricted ourselves to finite dimensional Lie algebras. The cohomology of finite-dimensional reductive (in particular, semi-simple) Lie algebras over a field of characteristic 0 has been investigated completely. But, up to now, only a few general assertions are known about the cohomology of solvable Lie algebras. 

	In 1983, L. J. Santharoubane computed explicitly the cohomology of Heisenberg Lie algebras (see \cite{San83}). In 2005, H. Pouseele calculated the cohomology of one-dimensional central extension of Heisenberg Lie algebras including the diamond Lie algebras (see \cite{Pou05}). For some other results about Lie algebras cohomology we refer the reader to the references \cite{ACJ97} or\cite{PT09}.

For the quadratic Lie algebras the description of their cohomology could be done in particular ways by using the skew-symmetric derivations or the super-Poisson bracket (see \cite{MR85}, \cite{PU07}, \cite{DPU12} and \cite{Duo14}).

In this paper, we will be concerned about the cohomology of some families of Lie algebras and quadratic Lie algebras.  The paper will be organized as follows: Section \ref{Sec 1} is devoted to some preliminary concepts and results. Section \ref{Sec 2} will describe the cohomology of all $MD\left( n,1 \right)-$algebras ($n\ge 2$) which were completely classified in \cite{LHNCN16}. 
Finally, Section \ref{Sec 3} will compute the second Betti numbers of the generalized complex diamond Lie algebras by using the super-Poisson bracket. 

\section{Preliminaries}\label{Sec 1}

In this section, we recall some usual notations, basic concepts and well-known results  which will be used later. For details we refer the reader to \cite{Pos86}, \cite{LHNCN16}, \cite{DPU12} and \cite{MR85}.

\subsection*{Notations}

Throughout this paper, we will use the usual notations as follows.
	\begin{itemize}
		\item Let $k, n$ be two integer numbers, $0 \leq k \leq n$. The classical notation
$\binom{n}{k}$ means the number of combinations of size $k$ from $n$. When $k$ is negative, it is assumed in accordance with the general conventions that $\binom{n}{k}: = 0$.
		\item The notation $\mathbb{K}$ will mean a field of characteristic zero. In this paper, $\mathbb{K}$ is always $\mathbb{R}$ or $\mathbb{C}$.
		\item For a Lie algebra $\mathfrak{g}$, we always denote by $\mathfrak{g}^*$ the dual space of $\mathfrak{g}$. Unless otherwise stated, for any basis $\{X, Y, \ldots \}$ in $\mathfrak{g}$, by $\{X^*, Y^*, \ldots \}$ we mean the dual basis of the first one in $\mathfrak{g}^*$.
		\item Let $\mathfrak{g}$ be an $n$-dimensional Lie algebra over over a field $\mathbb{K}$. By ${\Lambda}^{k}_{\mathbb{K}} (\mathfrak{g}^*)$ we always mean the space of all skew-symmetric $k$-linear forms (with $\mathbb{K}$-values) on $\mathfrak{g}$\, ($0 \leq k \leq n$). The notation $\Lambda_{\mathbb{K}}(\mathfrak{g}^*)$ will mean the algebra of external forms on $\mathfrak{g}$. Moreover, $\wedge$ will mean the external product on $\Lambda_{\mathbb{K}}(\mathfrak{g}^*)$ and $\iota_{X}$ is the interior product on $\Lambda (\mathfrak{g}^*)$ associated to one fixed element $X \in \mathfrak{g}$.
	\end{itemize}	 

\subsection{Cohomology of Lie algebras}

Let $\mathfrak{g}$ be a finite-dimensional Lie algebra over a field $\mathbb{K}$ of characteristic zero, and let $V$ be a module over $\mathfrak{g}$ having, as a vector space over $\mathbb{K}$, a finite dimension (in other words, $V$ is a space  of some finite-dimensional representation $\rho$ of $\mathfrak{g}$).

\begin{defn}[\cite{Pos86}]
A function $\omega=\omega (X_0, X_1, \ldots, X_{k-1} )$ of $k$ independent variables 
$X_0, X_1, \ldots , X_{k-1} \in \mathfrak{g}\,\, (k \in \mathbb{N} \setminus \{0 \})$, which takes on values in the module $V$, is said a {\it  $k$-cochain} of the Lie algebra $\mathfrak{g}$ over $V$ if $\omega$ is skew-symmetric $k$-linear. All $k$-cochains constitute in an obvious way a vector space $C^k (\mathfrak{g},V)$.
\end{defn}

By definiton, $C^1(\mathfrak{g}, V) = Hom_{\mathbb{K}}(\mathfrak{g}, V)$. It is assumed, in accordance with the general conventions, that $C^0(\mathfrak{g}, V) = V$.

For any cochain  $\omega \in C^k (\mathfrak{g},V)$ and $X_0, X_1, \ldots ,X_k \in \mathfrak{g}$, we set
\begin{multline*}
\delta \omega (X_0, X_1,\ldots ,X_k) = \sum\limits_{0\le i\le k} (-1)^i \rho (X_i) \omega (X_0, X_1, \ldots , \widehat{X_i} ,\ldots , X_k) \,+
  \\ 
+
 \sum\limits_{0\le i<j\le k} (-1)^{i+j} \omega ([ X_i, X_j], X_0, X_1, \ldots ,\widehat{X_i},\ldots , \widehat{X_j}, \ldots , X_k),
\end{multline*}
where the symbol  $\text{ } \widehat{ } \text{ }$  over the independent variable means that the latter must be omitted.
It is clear that $\delta \omega$ is a $(k+1)$-cochain for any $\omega \in C^k (\mathfrak{g},V)$. 

\begin{defn}[\cite{Pos86}]
The mapping 
 $$\delta : C^k (\mathfrak{g},V) \to C^{k+1}(\mathfrak{g},V); \omega \mapsto \delta \omega$$ 
 is obviously linear and it is called the {\it ($k$-th) coboundary operator}.
\end{defn}

The basic property of the mapping $\delta$ is that if twice repeated it is zero: $\delta \circ  \delta = 0$. 
A cochain $\omega$ for which $\delta \omega = 0$ is called a {\it cocycle} and a cochain $\omega$ of the form
$\delta \varphi$ is a {\it coboundary}. All  cocycles form a subspace $Z^k (\mathfrak{g},V)$ of $C^k (\mathfrak{g},V)$. All  coboundaries (for $k>0$) also form a subspace $B^k (\mathfrak{g},V)$ of $C^k (\mathfrak{g},V)$.
The relation $\delta \circ  \delta = 0$ means that
$$B^k (\mathfrak{g},V) \subset Z^k (\mathfrak{g},V); k > 0.$$
So that we have the following definition.

\begin{defn}[\cite{Pos86}]
The vector space
$$H^k (\mathfrak{g},V)= Z^k (\mathfrak{g},V)/ B^k (\mathfrak{g},V)$$
is called {\it $k-$th cohomology space} (or {\it $k-$th cohomology}, for brevity) of $\mathfrak{g}$ {\it with coefficients in 
$V$}. For any cocycle $\omega \in Z^k (\mathfrak{g},V)$, we will denote by $[\omega]$ the cohomology class of $\omega$ in 
$H^k (\mathfrak{g},V)$.
\end{defn}

When $k=0$ we agree to assume that
$H^0 (\mathfrak{g},V)= Z^0 (\mathfrak{g},V)$
so that $H^0 (\mathfrak{g},V)$ is nothing than a subspace of the module $V$ consisting of 
the {\it invariant  elements} of $V$, i.e. of elements $v \in V$ 
such that $\rho (X) (v) = 0$ for any $X \in \mathfrak{g}$.

In fact, computing cohomology of Lie algebras, one often considers some special cases of vector space  $V$ and representation $\rho $ of $\mathfrak{g}$ in $V$. In particular, when $V=\mathfrak{g}$ and $\rho = ad$  is the adjoint representation of  $\mathfrak{g}$, we get the assertion that $H^1 (\mathfrak{g},\mathfrak{g})$ is commonly used to describe the space of the outer derivations of $\mathfrak{g}$.




The special case $V = \mathbb{K}$ is one of the most remarkable cases of cohomology of Lie algebras and 
$H^k (\mathfrak{g},\mathbb{K})$ said to be ($k$-th) cohomology space of $\mathfrak{g}$ 
{\it with trivial coefficients}. In this case, 
$C^0 (\mathfrak{g},\mathbb{K})= \mathbb{K}$ and 
$C^k (\mathfrak{g},\mathbb{K}) = {\Lambda}^{k}_{\mathbb{K}} ( \mathfrak{g}^* )$. Moreover, the coboundary operator $\delta$ is given as follows:

$\delta _k f (X_0, \ldots , X_k) =$

\hskip 1.5cm $= \sum\limits_{ 0 \leq i<j \leq k} (-1)^{i+j} f (\left[ X_i, X_j \right], X_0, \ldots ,\widehat{X}_i,\ldots , \widehat{X}_j, \ldots , X_k ).$

\begin{defn}[\cite{MR85}]\label{defn 1.1.4}
The dimensional number of the cohomology $H^k (\mathfrak{g}, \mathbb{K})$ 
of the Lie algebra $\mathfrak{g}$ is called {\it the $k$-th Betti number} of $\mathfrak{g}$ and denoted by 
$b_k(\mathfrak{g})$.
\end{defn}

Note that vector spaces are isomorphic if they have the same dimension. Therefore, the description of the cohomology with trivial coefficients (in $\mathbb{K}$) of Lie algebras, in essence, is just calculating their Betti numbers. Moreover, the dimension of the cohomology can be calculated via the dimension of $(\mathfrak{g},V)$. Namely, we have the following result.

\begin{prop}[see, for example, \cite{MR85}] \label{prop1.1.5}
For any Lie algebra $\mathfrak{g}$ and $\mathfrak{g}$-module $V$, we have
\begin{enumerate}[(i)]
	\item $ \dim H^k (\mathfrak{g},V) = \dim Z^k (\mathfrak{g},V) + \dim Z^{k-1} (\mathfrak{g},V)-
\binom{n}{k-1} \dim V$.
	\item $ \dim H^k (\mathfrak{g},V) = \binom{n}{k} \dim V- \dim B^k (\mathfrak{g},V)- \dim B^{k+1} (\mathfrak{g},V). \hfill \square $
\end{enumerate}
\end{prop}

When $V = \mathbb{K}$, as an immediate consequence of Proposition \ref{prop1.1.5} we get following corollary which is very useful for us later.

\begin{cor}\label{cor1.1.6}
For any Lie algebra $\mathfrak{g}$, we have
\begin{enumerate}[(i)]
	\item  $ b_k (\mathfrak{g}) = \dim Z^k (\mathfrak{g},\mathbb{K}) 
	+ \dim Z^{k-1} (\mathfrak{g},\mathbb{K})-\binom{n}{k-1}$.
	\item  $ b_k (\mathfrak{g}) = \binom{n}{k}
- \dim B^k (\mathfrak{g},\mathbb{K})- \dim B^{k+1} (\mathfrak{g},\mathbb{K}). \hfill \square $
\end{enumerate}
\end{cor}
It is to be noticed that for the decomposable Lie algebras, i.e. the Lie algebras of the form $\mathfrak{g} \oplus \mathfrak{h}$, the computation of their cohomology can be reduced to that of $\mathfrak{g}$ and $\mathfrak{h}$ by using the K{\mbox{$\ddot{u}$}}nneth formula (see \cite{Sol02}). Namely, we have the following proposition.

\begin{prop}[\cite{Sol02}, Proposition 2.8]\label{prop1.1.7}
Let $\mathfrak{g}, \mathfrak{h}$ be finite-dimensional Lie algebras over a field $\mathbb{K}$. Then

\centerline{$H^* (\mathfrak{g}\oplus \mathfrak{h},\mathbb{K}) \cong
H^* (\mathfrak{g},\mathbb{K}) \otimes _{\mathbb{K}} H^*  (\mathfrak{h},\mathbb{K}).$}  \hfill $\square$
\end{prop}




 


\subsection{$MD(n, m)$-algebras and $MD(n, m)$-class}

\begin{defn}[\cite{LHNCN16}]\label{defn 1.2.1}
Let $n, m$ be natural numbers ($0 \leq m < n$) and $\mathbb{K} = \mathbb{R}$. Assume that $\mathfrak{g}$ is an $n$-dimensional solvable real Lie algebra such that  the first derived ideal $\mathfrak{g}^1 := [\mathfrak{g}, \mathfrak{g}]$ is $m$-dimensional one. Denote by $G$ the connected, simply connected Lie group associated to $\mathfrak{g}$. Then, $\mathfrak{g}$ is called an {\it $MD( n,m)$-algebra} if all coadjoint orbits of $G$ (i. e. the orbits in the coadjoint representation of $G$) are the ones of zero or maximal dimension. The class of all $MD(n, m)$-algebras is called the {\it $MD(n, m)$-class}.
\end{defn}

\subsection{Quadratic Lie algebras}

\begin{defn}[\cite{MR85}]\label{defn 1.3.1}
Let $\mathbb{K}=\mathbb{C}$ and $\mathfrak{g}$ be a complex Lie algebra. Then $\mathfrak{g}$ is called {\it a quadratic Lie algebra} if it is endowed with a non-degenerate invariant symmetric bilinear form (with complex values) $B$, i.e. $B$ satisfies the following conditions:
\begin{enumerate} [(i)]
\item $B(X,Y)=B(Y,X)$ for all $X,\ Y\in \mathfrak{g}$,
\item if $B(X,\mathfrak{g}) =0$ then $X=0$,
\item $B([X,Y],Z)=B(X,[Y,Z])$ for all $X,\ Y,\ Z\in \mathfrak{g}$.
\end{enumerate}
\end{defn}

\begin{defn}[\cite{PU07}]\label{defn 1.3.2}
Let $(\mathfrak{g},B)$ be a quadratic Lie algebra over $\mathbb{C}$. We choose in $\mathfrak{g}$ one basis $\{ X_1, \ldots , X_n\} $ with the dual basis $\{X^* _1, \ldots , X^*_n \}$. Assume that $\{ Y_1, \ldots ,Y_n\} $ is the basis of $\mathfrak{g}$ which is given by 
$B(Y_i, .) = X^*_i,\ 1 \leq i \leq n$.
Then, {\it the super-Poison bracket} on $\Lambda_{\mathbb{C}} ( \mathfrak{g}^*)$ is defined by
$$\{\Omega ,\Omega  '\}= (-1)^{k+1} \sum\limits_{1 \leq i, j \leq n} B (Y_i, Y_j) \iota _{X_i} (\Omega )\wedge \iota  _{X_j} (\Omega ')$$
for all \,$\Omega \in {\Lambda}^{k}_{\mathbb{C}} ( \mathfrak{g}^* ), \, \, \Omega ' \in \Lambda_{\mathbb{C}} (\mathfrak{g}^*)$.
\end{defn}

\begin{defn}[see \cite{PU07} or \cite{DPU12}]\label{defn 1.3.3}
The 3-form $I$ on $\mathfrak{g}$ which is defined by 
\[I(X,Y,Z)=B([X,Y],Z); \forall X,\ Y,\ Z \in \mathfrak{g}.\]
is called {\it the 3-form associated to $(\mathfrak{g}, B)$}.
\end{defn}

\begin{prop}[\cite{DPU12}]\label{prop1.3.4}
Via the 3-form $I$ associated to $(\mathfrak{g}, B)$, the coboundary operator $\delta $ of $(\mathfrak{g},B)$ is exactly represented as follows $-\{I, .\}$. That means \,\,
$\delta \Omega  =  - \{ I,\Omega \} $, for all  $\Omega \in \Lambda ( \mathfrak{g}^* )$. \hfill $\square$
\end{prop}

\section{Cohomology of $MD(n, 1)$-algebras}\label{Sec 2}
	
In this section, we will describe explicitly the cohomology of $MD\left( n,1 \right)$-class ($n \geq 2$). This class was completely classified in \cite{LHNCN16}. We emphasize that, throughout  this section, the basic field $\mathbb{K}$ is always $\mathbb{R}$ and we only consider the cohomology with real (trivial) coefficients. Before coming to the computations, we introduce this classification.

\subsection{The classification of all $MD(n, 1)$-algebras}

The $MD\left( n,1 \right)$-class includes only the Lie algebra of the group of affine transformations of the real straight line, the real Heisenberg Lie algebras and their trivial extensions by the real commutative Lie algebras. Namely, we have the following classification.

\begin{prop} [\cite{LHNCN16}]\label{ClassMD(n,1)}
		$MD(n,1)$-class coincides with the class of all real solvable Lie algebras whose the first derived ideal is 
		$1-$dimensional, moreover $MD(n,1)$ includes only the Lie algebra of the group of affine transformations of the real straight line, the real Heisenberg Lie algebras and their direct extensions by the real commutative Lie algebras. In other words, if $\mathfrak{g}$ is a $n$-dimensional real solvable Lie algebra whose the first ideal $\mathfrak{g}^1:= [\mathfrak{g}, \mathfrak{g}]$ is $1$-dimensional $(2 \leq n \in \mathbb{N})$ then $\mathfrak{g}$ is an $MD(n,1)$-algebra and $\mathfrak{g}$ is isomorphic to one and only one of the following Lie algebras.
		\begin{itemize}
		\item[(i)] The Lie algebra ${\rm aff} (\mathbb{R})$ of the group of all affine transformations on $\mathbb{R}$; $n = 2$.
		\item[(ii)] ${\rm aff} (\mathbb{R}) \oplus {\mathbb{R}}^{n-2}$; $n > 2.$
		\item[(iii)] The real Heisenberg Lie algebra $\mathfrak{h}_{2m + 1}$;\, $3 \leq n = 2m + 1.$
		\item[(iv)] $\mathfrak{h}_{2m+1} \oplus {\mathbb{R}}^{n-2m-1}$; $3 \leq 2m + 1 < n.$ \hfill $\square$
		\end{itemize}
\end{prop}

\subsection{Cohomology of all $MD(n, 1)$-algebras} 
First of all, we recall the definition of the Lie algebra ${\rm aff} (\mathbb{R})$ and the real Heisenberg Lie algebra $\mathfrak{h}_{2m + 1}$\, ($3 \leq n = 2m + 1)$. Namely, we have
\begin{itemize}
	\item  ${\rm aff}(\mathbb{R}): = span \left\{ X,Y \right\}$ with $[X,Y] = Y$.
	\item  $\mathfrak{h}_{2m+1}:= span \{ Z, X_1, \ldots, X_m, X_{m + 1}, \ldots, X_{2m} \},
[ X_i, X_{m + i} ] = Z$ for all $i = 1, 2, ..., m$. The other Lie brackets are trivial. 
\end{itemize}

The following theorem give the Betti number of all $MD(n, 1)$-algebras. In fact, the cohomology of the Heisenberg Lie algebras was described by L. J. Santharoubane in 1983 (see \cite{San83}). However, for the sake of completeness,
we will reintroduce here without the proof.
 
\begin{thm}\label{Thm2.2.1}
	Notations being as above, the Betti numbers of all $MD(n, 1)$-algebras are described as follows.
\begin{enumerate} [(i)] 
	\item		
	$b_1({\rm aff}(\mathbb{R})) = 1$, \, $b_2(\rm{aff}(\mathbb{R})) = 0$.
	\item
${b_k}\left({\rm aff}(\mathbb{R}) \oplus {\mathbb{R}^{n - 2}}\right) = \binom{n-1}{k}; \ 0 \leq k \leq n.$
	\item

$b_k \left(\mathfrak{h}_{2m+1}\right) = b_{n - k} \left(\mathfrak{h}_{2m+1}\right) = \binom{2m}{k} - \binom{2m}{k-2}; \, 0 \leq k \leq m.$
	\item
$b_k \left( \mathfrak{h}_{3}\oplus  \mathbb{R}^{n - 3}\right) =  \binom{n-1}{k} + \binom{n-2}{k-2}; \, n>3, 0 \leq k \leq n.$
	
$b_k \left(\mathfrak{h}_{2m+1}\oplus  \mathbb{R}^{n - 2m -1}\right) = b_{n-k} \left(\mathfrak{h}_{2m+1}\oplus  \mathbb{R}^{n - 2m -1}\right) =$

\hskip 1cm $= \binom{n-1}{k} - \binom{n-1}{k-2}; n > 2m + 1 > 3, \, \, 0 \leq k \leq m$. 

$b_k \left(\mathfrak{h}_{2m+1}\oplus  \mathbb{R}^{n - 2m -1}\right) = b_{n-k} \left(\mathfrak{h}_{2m+1}\oplus  \mathbb{R}^{n - 2m -1}\right) =$ 

\hskip 1.5cm $=\sum\limits_{i=0}^{\min \{k,2m+1\}} \left( \binom{2m}{i - \left[\frac{i}{m+1}\right]} - \binom{2m}{i + 3 \left[\frac{i}{m+1}\right] - 2} \right)  \binom{n -2m - 1}{k-i};$\\ 
$n > 2m + 1 > 3, \, \, m+1 \leq k \leq \left[ \frac{n}{2} \right]$. 
\end{enumerate}
\end{thm}

\noindent {\bf The Proof of Theorem 2.2.1}

The proof of Part {\it (i)} is straight-forward by an elementary computation. Moreover, Part {\it (iii)} has proved in \cite{San83}. So we need only prove Parts {\it (ii), (iv)} and {\it (v)}.
	\begin{itemize}	
		\item {\it Proof (ii)}: For ${\rm aff} (\mathbb{R}) \oplus {\mathbb{R}}^{n-2}\, (n > 2)$. 
		It is clear that this algebra accepts the basis $\left\{ X, Y, Z_1, Z_2,\ldots , Z_{n - 2} \right\}$ 
with only one non-trivial bracket $[X, Y] = Y$. 

Direct computation by using the representation the coboundary operator shows that
		
$\delta (Z_{i_1}^* \wedge Z_{i_2}^*  \wedge \ldots \wedge Z_{i_k}^*) = 0,\, \delta ( X^* \wedge Z_{i_1}^* \wedge \ldots \wedge Z_{i_{k - 1}}^*) = 0$,

$\delta \left( X^* \wedge Y^* \wedge Z_{i_1}^* \wedge \ldots \wedge Z_{i_{k - 2}}^* \right) = 0$,

$\delta \left( Y^* \wedge Z_{i_1}^* \wedge \ldots \wedge Z_{i_{k - 1}}^* \right) =  - {X^*} \wedge Y^* \wedge Z_{i_1}^*  \wedge \ldots \wedge Z_{i_{k - 1}}^*$,		

where $1 \leq i_1  < \ldots  < i_{k} \leq n - 2$.		

Therefore, we obtain

$B^k \left({\rm aff}(\mathbb{R}) \oplus \mathbb{R}^{n - 2}\right) =$

\hskip 0.5cm $=span\left\{ X^* \wedge Y^* \wedge Z_{i_1}^* \wedge \ldots  \wedge Z_{i_{k - 2}}^*; \right.$ 

 \hskip 5.5cm $\left. 1 \leq i_1  < \ldots  < i_{k - 2} \leq n - 2 \right\}.$

$Z^k \left({\rm aff}(\mathbb{R}) \oplus \mathbb{R}^{n - 2}\right) =$

\hskip 0.5cm $= span \left \{ Z_{i_1}^* \wedge \ldots \wedge Z_{i_k}^*,  
  X^* \wedge Z_{i_1}^*  \wedge \dots \wedge Z_{i_{k - 1}}^*,\right.$
 
 \hskip 1.5cm $\left.  X^* \wedge Y^*  \wedge Z_{i_1}^* \wedge \ldots \wedge Z_{i_{k - 2}}^*; 
1 \leq i_1  < \ldots < i_k \leq n-2 \right \}.$

$H^k \left({\rm aff}(\mathbb{R}) \oplus \mathbb{R}^{n - 2}\right) =$

\hskip 0.5cm $= span\left\{ \left[ Z_{i_1}^* \wedge \ldots \wedge Z_{i_k}^* \right], 
  \left[ X^* \wedge Z_{i_1}^* \wedge \ldots \wedge Z_{i_{k - 1}}^* \right]; \right. $

 \hskip 5.5cm $ \left.  1 \leq  i_1 < ... < i_k  \leq n-2 \right\}.$
 
Therefore, we get

\hskip 0.5cm $ b_k \left({\rm aff}(\mathbb{R}) \oplus \mathbb{R}^{n - 2}\right) = 
\binom{n-2}{k} + \binom{n-2}{k-1} = \binom{n-1}{k}; 0 \leq k \leq n.$

 	\item {\it Proof (iv)}: Note that if $\mathfrak{g}$ is the Heisenberg algebra or its trivial extention then $\mathfrak{g}$ and, by the Poincare duality, $ H^k ( \mathfrak{g}) \cong H^{n-k} (\mathfrak{g})$, where $n = \dim \mathfrak{g}$. That means $b_k(\mathfrak{g}) = b_{n-k}(\mathfrak{g})$. Therefore, it needs only to study 
$H^k (\mathfrak{g})$ with $k \leq \left[\dfrac{n}{2}\right]$.

By virtue of Propositon \ref{prop1.1.7}, we get 

$H^k (\mathfrak{h}_{2m+1}\oplus \mathbb{R}^{n-2m-1}, \mathbb{R}) = $
 
\hskip 0.5cm $=\underset{i=0}{\overset{k}{\mathop{\oplus }}} \,
\left( H^i ( \mathfrak{h}_{2m+1}, \mathbb{R}) \otimes _{\mathbb{R}} H^{k-i} (\mathbb{R}^{n-2m-1}, \mathbb{R}) \right).$
\vskip 0.5cm
Therefore, we obtain

$b_k ( \mathfrak{h}_{2m+1} \oplus \mathbb{R}^{n-2m-1}) = \sum\limits_{i=0}^{k} {b_i( \mathfrak{h}_{2m+1}). b_{k-i} ( \mathbb{R}^{n-2m-1})}.$

 \begin{enumerate}[a)]
\item 
 For $m=1, 3 \leq k \leq n-1$, we see at once that  

$b_k (\mathfrak{h}_{3} \oplus \mathbb{R}^{n-3})
= \binom{n-3}{k}+2 \binom{n-3}{k-1}+2\binom{n-3}{k-2}+\binom{n-3}{k-3}.$ 

It follows immediately that
$b_k (\mathfrak{h}_{3} \oplus \mathbb{R}^{n-3}) = \binom{n-1}{k}+ \binom{n-2}{k-2}.$

It is easy to check that this formula is accurate for $k=0,1,2.$
\item 
 For $m>1,0 \leq k \leq m$,  we get

$ b_k (\mathfrak{h}_{2m+1} \oplus \mathbb{R}^{n - 2m -1})= \binom{n-1}{k} - \binom{n-1}{k-2}.$ 

\item
For $m+1 \leq k \leq \left[ \dfrac{n}{2} \right]$, we obtain

$b_k ( \mathfrak{h}_{2m+1} \oplus \mathbb{R}^{n-2m-1}) = $

\hskip 0.7cm $ = \sum\limits_{i=0}^{m}{\left( \binom{2m}{i} - \binom{2m}{i-1}\right)\binom{n-2m-1}{k-i}} + $

\hskip 1.5cm $ +\sum\limits_{i=m+1}^{\min \{k,2m+1\}} \left( \binom{2m}{2m+1-i} - \binom{2m}{2m+1-i-2} \right) \binom{n-2m-1}{k-i}$

\vskip0.4cm

Direct computation shows that 

$b_k ( \mathfrak{h}_{2m+1} \oplus \mathbb{R}^{n-2m-1}) = $
 
\hskip 1cm $ = \sum\limits_{i=0}^{\min \{k,2m+1\}}\left( \binom{2m}{i-\left[\frac{i}{m+1}\right]} - \binom{2m}{i+3\left[\frac{i}{m+1}\right]-2} \right) \binom{n-2m-1}{k-i};$\\
$n > 2m + 1 > 3, \, \, m+1 \leq k \leq \left[ \frac{n}{2} \right]$.

\end{enumerate}

\end{itemize}
The proof of Theorem \ref{Thm2.2.1} is complete. \hfill $\square$

\section{Cohomology of a special family \\ of quadratic Lie algebras}\label{Sec 3}
In this section, we will apply the super-Poisson bracket (see \cite{PU07} for more details) to compute the second cohomology of a special family of quadratic Lie algebras which is called the family of generalized diamond Lie algebras. 
\subsection{The generalized diamond Lie algebras}

Recall that the diamond Lie algebra $\mathfrak{D}_{4}$ is the 4-dimensional one which is given as follows:
$$\mathfrak{D}_{4}:= span \{ X_0, X_1, Y_0, Y_1 \}; [Y_0,X_1]=X_1, [Y_0,Y_1]=-Y_1, [X_1,Y_1]=X_0$$ 
and the other Lie brackets are trivial.

Generalizing this algebra, for any 
$\Lambda: = (\lambda_1, \lambda_2, \ldots , \lambda_n) \in \mathbb{K}^n$ ($0 < n \in \mathbb{N}$), we consider the $(2n + 2)$-dimensional Lie algebras $\mathfrak{D}_{2n+2}(\Lambda) \,$ over $\mathbb{K}$ which is defined as follows:
$$\mathfrak{D}_{2n+2}(\Lambda) := span \{ X_0,\ldots , X_n, Y_0, \ldots , Y_n \}$$
with non-trivial Lie brackets are given by
$$[ Y_0, X_i] =  \lambda _i X_i, [Y_0, Y_i] = - \lambda _i Y_i, [X_i, Y_i] = \lambda _i X_0;  1\le i\le n.$$ 

We obtain one infinite family of the $(2n + 2)-$dimensional Lie algebras 
$$\{ \mathfrak{D}_{2n+2}(\Lambda); \, \Lambda \in \mathbb{K}^n \}$$ 
which is numbered by $\Lambda \in \mathbb{K}^n$. Note that when $n=1, \Lambda \equiv \lambda_1 \neq 0$ then $\mathfrak{D}_{2+2}(\lambda _1)$ is exactly isomorphic to the diamond Lie algebra $\mathfrak{D}_4$. 

Moreover, $\mathfrak{D}_{2n+2}(\Lambda)$ can be endowed with a non-degenerate invariant symmetric bilinear $2$-form $B$ which is defined as follows:
$$B( X_i, Y_i) = 1,\ B(X_i, X_j) = B(X_i, Y_j) = B(Y_i, Y_j) = 0; 0 \le i \neq j \le n.$$
It is easy to verify that $(\mathfrak{D}_{2n+2}(\Lambda), B)$ is a quadratic Lie algebra over $\mathbb{K}$.
Motivated by this fact, we introduce the following notion.

\begin{defn}
The quadratic Lie algebra $(\mathfrak{D}_{2n+2}(\Lambda), B)$ is called {\it the generalized diamond Lie algebra} (over $\mathbb{K}$). When  $\mathbb{K}= \mathbb{R}$ (or $\mathbb{K}= \mathbb{C}$, respectively) then  it is called {\it the generalized real (or complex, respectively) diamond Lie algebra.}
\end{defn}

We emphasize that, the generalized real diamond Lie algebras ($\mathbb{K} = \mathbb{R}$) were mentioned by A. Medina and P. Revoy (see \cite{MR85}) in 1985 with the restricted condition that $\Lambda \in \mathbb{R}_{+}^n$, i.e. $0 < \lambda _i \in \mathbb{R},\ i = 1, 2, \ldots, n$. In fact, the authors have calculated the second cohomology with real (trivial) coefficients of $(\mathfrak{D}_{2n+2}(\Lambda), B)$ for $\Lambda \in \mathbb{R}_{+}^n$ via the description of the skew-symmetric derivations. 

In the next subsection, we will compute the second Betti number of the generalized complex diamond Lie algebras ($\mathbb{K} = \mathbb{C}$) with complex (trivial) coefficients for $\Lambda \in \mathbb{C}^n$ by using the super-Poison bracket. We emphasize that the method used here is not only new, but also more effective than the method of A. Medina and P. Revoy in \cite{MR85} because it makes the computations become simpler and briefer.

\begin{rem}\label{rem3.1.2} We have some remarks as follows.
\begin{enumerate}
	\item[(i)] Let $\sigma$ be a permutation of the set $\{1, 2, \ldots, n\}$. We set 
	$$\sigma (\Lambda) : = (\lambda_{\sigma (1)}, \ldots \lambda_{\sigma (n)}).$$
Assume that $\tilde{\sigma}$ is the map induced by $\sigma$ as follows
$$\tilde{\sigma} : \mathfrak{D}_{2n + 2} (\Lambda) \longrightarrow \mathfrak{D}_{2n + 2} (\sigma(\Lambda));\, 
\tilde{\sigma}(X_i):=X_{\sigma (i)}, \tilde{\sigma}(Y_i):=Y_{\sigma (i)},$$
where $1 \leq i \leq n$. It is obvious that $\tilde{\sigma}$ is an Lie algebra isomorphism. Moreover, $B$ is invariant for $\tilde{\sigma}$, then $\tilde{\sigma}$ is also a quadratic Lie algebra isomorphism.
	\item[(ii)] Now we consider some $\Lambda = (\lambda_1, \dots, \lambda_n) \in \mathbb{K}^n$. Assume that 
$\lambda _i =0$ for $i \in I \subset \{1, 2, \ldots, n\}$. In view of Part (i), we can assume, without loss of generality, that 
$$\Lambda = (\Lambda_m, 0, 0, \ldots ,0),$$
where  
$\Lambda _m = (\lambda_1, \ldots , \lambda_m)  \in \left(\mathbb{K}  \backslash \{ 0\}\right)^m; \, 0 \leq m \leq n$.
It is easy to check that $(\mathfrak{D}_{2n + 2}(\Lambda), B)$ is decomposable. Namely, we have  
$$(\mathfrak{D}_{2n+2}(\Lambda), B)=(\mathfrak{D}_{2m+2}(\Lambda _m), B_m) \oplus \mathbb{K}^{2(n-m)},$$
where $B_m$ is the restriction of $B$ on $\mathfrak{D}_{2m+2} \subset \mathfrak{D}_{2n+2}$. Therefore, 
the computation of cohomology of $(\mathfrak{D}_{2n+2}(\Lambda), B)$ can be reduced to that of $(\mathfrak{D}_{2m+2}(\Lambda _m), B_m)$ by applying the K{\mbox{$\ddot{u}$}}nneth formula in Proposition \ref{prop1.1.7}. 
\end{enumerate}
\end{rem}

\subsection{The second cohomology of the generalized complex diamond Lie algebras}	
Consider some algebra from the family of 
$(\mathfrak{D}_{2n+2}(\Lambda), B)$ over $\mathbb{C}$ in which one fixed $\Lambda = (\lambda_1, \lambda_2, \ldots , \lambda_n) \in \mathbb{C} ^n$ has been selected. For the sake of convenience, will write $\mathfrak{D}_{2n+2}$ for $(\mathfrak{D}_{2n+2}(\Lambda), B)$ when no confusion can arise. 

Taking Remark \ref{rem3.1.2} into account, we need consider the case 
$\Lambda = (\lambda_1, \lambda_2, \ldots , \lambda_n) \in \left(\mathbb{C}  \backslash \{ 0\}\right)^n$, i.e. $\lambda_i  \neq 0$ for all $i = 1, 2, \dots, n$. 

Note that, in the given $\Lambda$, the non-zero complex numbers $\lambda_1, \lambda_2, ..., \lambda_n$ are not necessarily different. Now, we define the subset $\{a_1,\ a_2, \ \ldots ,\ a_k\}$ of $\mathbb{C}\backslash \{ 0\}$ by the following conditions:
\begin{enumerate}
	\item[(i)] $\forall \ i \in \{1, \ldots, n\}, \, \exists !\,j \in \{1, \ldots k\}$, such that $a_j = \lambda_i$ or $a_j = - \lambda_i$.
	\item[(ii)] $\forall \ j \in \{1, \ldots, k\}, \, \exists \,i \in \{1, \ldots n\}$, such that $a_j = \lambda_i$ or $a_j = - \lambda_i$.
\end{enumerate}

For each $i \in \{1, \ldots, k \}$, we denote by $p_i$ (respectively, $q_i$) the number of times that the value $a_i$ (respectively, $-a_i$) is repeated in  $\Lambda = (\lambda _1,\ \lambda _2, \ \ldots ,\ \lambda _n)$. Now we set $n_i : = p_i + q_i$ for every $i = 1, \dots k$. Note that $\sum\limits_{i = 1}^k {n_i}  = n$. 

It is clear that, for any given 
$\Lambda = (\lambda_1, \lambda_2, ..., \lambda_n) \in \left(\mathbb{C}  \backslash \{ 0\}\right)^n$, 
the numbers $n_i \, \ldots, n_k$ are well-defined and it is not difficulty to calculate them.
\vskip 0.5cm

The following theorem is the main result of this subsection.

\begin{thm}\label{Thm3.2.1}
	Notations being as above, if
$\Lambda = (\lambda_1, \lambda_2, \ldots , \lambda_n) \in \left(\mathbb{C}  \backslash \{ 0\}\right)^n$,	
	then the second Betti number of $(\mathfrak{D}_{2n+2}(\Lambda), B)$ is given by the following formula:
$$b_2 (\mathfrak{D}_{2n+2}) = \sum\limits_{i=1}^{k}{n_{i}^{2}} -1.$$
\end{thm}

\noindent {\bf The Proof of Theorem 3.2.1}

In this proof, for simplicity of notation, we will denote by
$$\{\alpha , \alpha _1, \ldots , \alpha _n, \beta , \beta _1, \ldots , \beta _n \}$$ 
the dual basis in $\mathfrak{D}_{2n+2}^*$ of 
$\{ X_0, \ldots, X_n, Y_0, \ldots, Y_n \}$ in $\mathfrak{D}_{2n+2}$. Moreover, we set
 $$V:=span\{ \alpha _i : 1\le i\le n \},  W :=span\{ \beta _i : 1\le i\le n\}.$$ 
It can easily be checked that the 3-form $I$ associated to $ \mathfrak{D}_{2n+2}$ is given by:  
$$I=\beta \wedge \sum\limits_{i=1}^{n} \lambda _i \alpha _i \wedge \beta _i.$$
By setting
$ \Omega _n := \sum\limits_{i=1}^{n} \lambda _i  \alpha _i \wedge \beta _i$, we have 

$B^2 (\mathfrak{D}_{2n+2}) = \, \{ \iota _X (I) : X\in \mathfrak{g}_{2n+2} \}$

\hskip 1.95cm $ = \, span\left\{ \beta \wedge \alpha _i, \beta \wedge \beta _i, \Omega _n :1\le i\le n \right\}.$

\begin{itemize}
	\item  For $n = 1$, by an easy computation we get $H^2 (\mathfrak{D}_4) = \{0\}$.
	\item  For $n > 1$, direct computation by using Proposition \ref{prop1.3.4} shows that	
\begin{enumerate} [(i)] 
\item	
$\{ I,\alpha  \wedge {\alpha _i}\}  = {\alpha _i} \wedge {\Omega _n} - {\lambda _i}\alpha  \wedge \beta  \wedge {\alpha _i}$,
\item	
 $\{ I,\alpha  \wedge {\beta _i}\}  = {\beta _i} \wedge {\Omega _n} + {\lambda _i} \alpha  \wedge \beta  \wedge {\beta _i}$,
\item	
	$\{ I,\alpha  \wedge \beta \}  = I$,
\item	
$\{ I,{\alpha _i} \wedge {\alpha _j}\}  = ({\lambda _i} + {\lambda _j})\beta  \wedge {\alpha _i} \wedge {\alpha _j}$,
\item	
 $\{ I,{\beta _i} \wedge {\beta _j}\}  =  - ({\lambda _i} + {\lambda _j})\beta  \wedge {\beta _i} \wedge {\beta _j}$,
\item	
 $\{ I,{\alpha _i} \wedge {\beta _j}\}  = ({\lambda _i} - {\lambda _j})\beta  \wedge {\alpha _i} \wedge {\beta _j}$.
\end{enumerate}	
\end{itemize} 

Therefore, $ Z^2 (\mathfrak{D}_{2n+2})$ is spanned by the union of the following sets:

\hskip 1cm $\left\{ {\beta  \wedge {\alpha _i}, \beta  \wedge {\beta _i},{\Omega _n}; \,1 \leq  i  < j \leq  n} \right\}$,

\hskip 1cm $\left\{ {\alpha _i} \wedge {\alpha _j}; \,{\lambda _i} + {\lambda _j} = 0, \,1 \leq  i  <  j  \leq  n \right\}$,

\hskip 1cm $\left\{ {\beta _i} \wedge {\beta _j}; \,{\lambda _i} + {\lambda _j} = 0, \, 1 \leq  i  <  j  \leq  n \right\}$,

\hskip 1cm $\left\{ {\alpha _i} \wedge {\beta _j}; \, {\lambda _i} - {\lambda _j} = 0, \,
 1\leq  i \leq n, 1\leq  j \leq n \right\}$.

It follows that the second cohomology $H^2 (\mathfrak{D}_{2n+2})$ is spanned by the union of the set of the cohomology classes in $Z^2 (\mathfrak{D}_{2n+2})$ modulo $B^2 (\mathfrak{D}_{2n+2})$ as follows:
 
\hskip 1cm $\left\{ \left[ {{\alpha _i} \wedge {\alpha _j}} \right]:{\lambda _i} + {\lambda _j} = 0, \,1 \leq  i < j \leq  n \right\}$,
 
\hskip 1cm $\left\{ \left[ {{\beta _i} \wedge {\beta _j}} \right]:{\lambda _i} + {\lambda _j} = 0, \, 1 \leq  i  < j \leq  n \right\}$,

\hskip 1cm $\left\{ \left[ {{\alpha _i} \wedge {\beta _j}} \right]:{\lambda _i} - {\lambda _j} = 0, \, 1 \leq i \leq n, 
\,1 \leq j \leq n \right\} \backslash \left\{ \left[ {\alpha _1} \wedge {\beta _1} \right]\right\}$.

Now, to compute the multiple number, we construct an relation $\sim$ on the set 
$\{(\lambda _1, 1), \, (\lambda _2, 2), \, \ldots, (\lambda _n, n)\}$ as follows: 
$$(\lambda _i, i) \sim (\lambda_j, j) \Longleftrightarrow \lambda_i = \pm \lambda_j.$$
It is very easy to check that \,$\sim$ \,is one equivalence relation. We obtain the quotient set as follows:
$$\{ N_1, N_2, \ldots , N_k \} := \{(\lambda _1, 1), \, (\lambda _2, 2), \, \ldots, (\lambda _n, n)\} \diagup \sim .$$ 
By renumbering, if necessary, we get that the cardinality of $N_i$ is exactly $n_i, \, i = 1, \ldots k$. Then, $p_i$ (respectively, $q_i$) is exactly the number of elements in $\Lambda = (\lambda_1, \, \lambda_2,  \ldots, \lambda_n)$ which are equal to $a_i$ (respectively, $- a _i$). Recall that $\sum\limits_{i = 1}^k {n_i}  = n$ and $p_i + q_i = n_i$. Therefore we obtain

 $b_2 (\mathfrak{D}_{2n+2}) = \dim Z^2 ( \mathfrak{D}_{2n+2}) - \dim B^2 ( \mathfrak{D}_{2n+2})$
 
\hskip 1.78cm $ = 2\sum\limits_{i=1}^k {p_i} {q_i} +  \sum\limits_{i=1}^k {p_i^2}  + \sum\limits_{i=1}^k {q_i^2}   - 1$
 
\hskip 1.78cm $= \sum\limits_{i=1}^k ({p_i} + {q_i})^2  - 1 = \sum\limits_{i=1}^k {n_i^2}  - 1$.
\vskip 0.5cm
\noindent The proof is complete. \hfill $\square$

\vskip 0.5cm
\noindent {\bf Concluding remark.}  We conclude the paper with some comments as follows:

	\begin{itemize} 
			\item Applying Theorem \ref{Thm3.2.1}, we can consider some particular cases. For example, if all $\lambda_i = \pm \,c$ for some non-zero constant $c \in \mathbb{C}$, then we obtain
		 $$b_2 ( \mathfrak{D}_{2n+2}) = n^2-1.$$
			\item Similarly, in the case when  
${\lambda _1},\ \ldots ,{\lambda _n}, - {\lambda _1}, \ \ldots , - {\lambda _n}$
 are different from each other then we get
$$ b_2 ( \mathfrak{D}_{2n+2}) =  n-1.$$
	\item We note that in 2014, the second author M. T. Duong \cite{Duo14} have computed explicitly not only $b_2(\mathfrak{D}_{2n + 2})$ but also $b_k (\mathfrak{D}_{2n+2})$ for all $k = 0, \ldots 2n+2$. But he made this in a restricted condition that $\lambda _i = 1, \, i = 1, 2, \dots n$.
	\item In the general case, the problem of computation of $b_k ( \mathfrak{D}_{2n+2})$ for all $k = 0, 1, 3, \ldots 2n+2$\,($n > 2$) is still open. We will consider this problem in the next paper. 
     \end{itemize}

\section*{Acknowledgements} The authors would like to thank the University of Economics and Law, VNU-HCMC and Ho Chi Minh City University of Education for financial supports. The authors would like to thank Professor Nguyen Van Sanh for his encouragement.

\end{document}